\documentclass[12pt, reqno]{amsart}
\usepackage{amsmath, amsthm, amscd, amsfonts, amssymb, graphicx, color}
\usepackage[bookmarksnumbered, plainpages, backref]{hyperref}

\theoremstyle{definition}

\theoremstyle{remark}

\numberwithin{equation}{section}

\begin{document}

\setcounter{page}{1}

\noindent\parbox{4.85in}{\hspace{0.1mm}\\[1.5cm]\noindent
\textcolor{green}{J. Nonlinear  Sci. Appl. 0 (0000), no. 0, 00--00}\\
$\frac{\rule{6.55in}{0.05in}}{{}}$\\
{\footnotesize \textcolor{red}{\textsc{\textbf{\large{T}}he
\textbf{\large{J}}ournal of \textbf{\large{N}}onlinear
\textbf{\large{S}}ciences \text{and}
\textbf{\large{A}}pplications}}\\
\textcolor[rgb]{0.00,0.00,0.84}{\textbf{http://www.tjnsa.com }}\\
$\frac{{}}{\rule{6.55in}{0.05in}}$}\\[.5in]}
\title[Some remarks on the action of Lusin area operator]{Some remarks on the action of Lusin area operator in Bergman spaces of the  unit ball}

\author[Romi Shamoyan, Haiying  Li]{Romi Shamoyan$^1$ and Haiying  Li$^2$$^{*}$}

\address{$^{1}$ Department of Mathematics, Bryansk State
University, Bryansk 241050, Russia.}

\email{\textcolor[rgb]{0.00,0.00,0.84}{rsham@mail.ru}}

\address{$^{2}$ College of Mathematics and
Information Science, Henan Normal University, Xinxiang 453007,
P.R.China}
\email{\textcolor[rgb]{0.00,0.00,0.84}{tslhy2001@yahoo.com.cn}}

\dedicatory{The research is financed by Asian Development Bank.
No. 2006-A171(Sponsoring information)}

\subjclass[2000]{Primary  ; Secondary .}


\date{Received: ; Revised: .
\newline \indent $^{*}$ Haiying  Li}

\begin{abstract}
We study the action of Lusin area operator on Bergman classes in
the unit ball, providing some direct generalizations of recent
results of Z. Wu.
\end{abstract} \maketitle

\section{Introduction }

Let as usual $\textmd{B}=\{z\in \mathbb{C}^{n}:\,|z|<1 \}$ be the
open unit ball of $\mathbb{C}^{n}$ and $\textmd{S}$ the unit
sphere plane of $\mathbb{C}^{n}$.
 Let  $d\upsilon$ be the normalized Lebesgue measure on  $\textmd{B}$ and $d \sigma$ the normalized rotation invariant
 Lebesgue measure on  $\textmd{S}$. We denote by $H(\textmd{B})$ as usual the class of all holomorphic functions on $\textmd{B}$.

 For any real parameter $\alpha$ we consider the weighted volume measure
  $$d \upsilon_{\alpha}(z)=(1-|z|^{2})^{\alpha}d \upsilon(z).$$

  Suppose $0<p<\infty$ and $\alpha>-1,$ the weighted Bergman space $A_{\alpha}^{p}$ consists of those functions $f \in H(\textmd{B})$
  for which $$\|f\|_{A_{\alpha}^{p}}^{p}=\int_{\textmd{B}}|f(z)|^{p}d \upsilon_{\alpha}(z)<\infty.$$

  Let $r>0$ and $z\in \textmd{B}$, the Bergman metric ball at $z$ is  defined as
  $$D(z,r)=\{w\in \textmd{B}:\, \beta(z,w)=\frac{1}{2}\log\frac{1+|\varphi_{z}(w)|}{1-|\varphi_{z}(w)|}<r\}.$$

  Where the involution $\varphi_{z}$ has the form
  $$\varphi_{z}(w)=\frac{z-P_{z}(w)-s_{z}Q_{z}(w)}{1-\langle w,z\rangle},$$
    where by $s_{z}=(1-|z|^{2})^{\frac{1}{2}},\,P_{z}$ is the orthogonal projection into the space spanned by $z\in \textmd{B}$, i.e.
  $P_{w}(z)=\frac{\langle w,z\rangle z}{|z|^{2}},\,P_{0}(w)=0$ and $Q_{z}=I-P_{z}$ (see \cite{8}). The volume of $D(z,r)$ is given by (see \cite{8})
  $$\upsilon(D(z,r))=\frac{R^{2n}(1-|z|^{2})^{n+1}}{(1-R^{2}|z|^{2})^{n+1}},$$
    where $R=\tanh(r).$ Set $|D(z,r)|=\upsilon(D(z,r)).$ For $w\in D(z,r),\,r>0,$ we have that (see, for example, \cite{8})

 $$(1-|z|^{2})^{n+1}\asymp (1-|w|^{2})^{n+1}\asymp |1-\langle z,w\rangle |^{n+1}\asymp |D(z,r)| \eqno(1) $$
 and
 $$|D(z,r)|^{\alpha+n+1}\asymp v_{\alpha}(D(z,r)). \eqno(2) $$

For any $\zeta \in \textmd{S}$ and $r>0$, the nonisotrpic metric
ball $Q_{r}(\zeta)$ is defined by (see \cite{8})

$$Q_{r}(\zeta)=\{z\in \textmd{B}:\, |1-\langle z,\zeta\rangle |^{\frac{1}{2}}<r\}.$$

A positive Borel measure $\mu$ on $\textmd{B}$ is called a
$\gamma-$Carleson measure if there exists a constant $C>0$ such that
$$\mu(Q_{r}(\zeta))\leq C r^{2\gamma} \eqno(3)$$
for all $\zeta\in \textmd{S}$ and $r>0$.

A well-known result about the $\gamma-$Carleson measure (see
\cite{7}) is that $\mu$  is  a $\gamma-$Carleson measure if and only
if $$\sup_{a\in
\textmd{B}}\int_{\textmd{B}}\bigg(\frac{1-|a|^{2}}{|1-\langle
z,a\rangle|^{2}}\bigg)^{\gamma}d
\mu(z)<\infty,\,\,\gamma>0.\eqno(4)$$

The area operator relates to the nontangential maximal function,
Littlewood-Paley operator, multipliers and tent space. It is very
useful in the harmonic analysis. On the unit disk, the boundedness
and compactness of the area operators was studied by Cohn and Wu
respectively on the Hardy space and the weighted Bergman space (see
\cite{1,7}).

Motivated by the results of \cite{1,7}, we define the area operator
on the unit ball as follows. Let $\mu$ be a positive measure on
$\textmd{B}$, we define Lusin area operator
$$G_{\mu,\sigma}(f)(\xi)=\int_{\Gamma_{\sigma}(\xi)}|f(z)|\frac{d\mu(z)}{(1-|z|)^{n}},\,\,\,f\in H(\textmd{B}).$$
Here $\Gamma_{\sigma}(\xi)$ is the corresponding Koranyi approach
region with vertex $\xi$ on $\textmd{S}$, i.e.
$$\Gamma_{\sigma}(\xi)=\{z\in \textmd{B}:\, |1-\langle z,\zeta\rangle|<\sigma(1-|z|^{2}),\,\sigma>1\}.$$
The purpose of this paper is to study the area operator on the
weighted Bergman space in the unit ball in $\mathbb{C}^{n}$.

All embedding theorems we prove in this paper for area operator in
the ball were obtained recently by Z. Wu for n=1 (case of unit
disk), as in \cite{7} our proofs are heavily based on various
properties of so-called sampling sequences or r-lattice $\{a_k\}$ in
the unit ball (see \cite{8}) and various estimates connected with
Bergman metric ball in $\mathbb{C}^n$  (\cite{8}). In next section
we collect preliminaries, in final section we provide formulations
and proofs of all main results of this paper.

 Throughout the paper,   constants are  denoted by $C$, they are  positive and may differ from one occurrence to the other.
 The notation $A\asymp B$ means that there is a positive constant $C$ such that $C^{-1}B\leq A\leq CB$ and
 $A\lesssim B$ if there is a positive constant $C$ such that $A<CB.$

\section{Preliminaries }

To state and prove our results, let's collect some nice properties
of the Bergman metric ball that will be used in this paper.

\vspace{2mm}

{\bf Lemma 1. (\cite{8})} There exists a positive integer $N$ such
that for any $0<r\leq 1$ we can find a sequence $\{a_{k}\}$ in
$\textmd{B}$ with the following properties:

\indent(1) $\textmd{B}=\bigcup_{k} D(a_{k},r)$;\\
 \indent(2) The sets  $D(a_{k},\frac{r}{4})$ are mutually disjoint;\\
\indent(3) Each point $z\in\textmd{B}$ belongs to at most $N$ of the
sets  $D(a_{k},2r)$. \vspace{2mm}

{\bf Remark 1.} If $\{a_{k}\}$ is a sequence from Lemma 1, according
to the result on page 78 of \cite{8}, there exist positive constants
$C_{1},C_{2}$ such that
$$C_{1}\int_{\textmd{B}}|f(z)|^{p}d \upsilon_{\alpha}(z)\leq\sum_{k=1}^{\infty}|f(a_{k})|^{p}(1-|a_{k}|^{2})^{n+1+\alpha}\leq C_{2}\int_{\textmd{B}}|f(z)|^{p}d \upsilon_{\alpha}(z)\eqno(5)$$

for all $f\in A_{\alpha}^{p}.$ Such a sequence will be called a
sampling sequence or r-lattice for $A_{\alpha}^{p}.$

\vspace{2mm}

{\bf Lemma 2. (\cite{8})} For every $r>0$ there exists a positive
constant $C_{r}$ such that
$$C_{r}^{-1}\leq\frac{1-|a|^{2}}{1-|z|^{2}}\leq C_{r},\,\,\,C_{r}^{-1}\leq \frac{1-|a|^{2}}{|1-\langle a,z\rangle|} \leq C_{r},$$
for all $a$ and $z$ in $\textmd{B}$ such that $\beta(a,z)<r.$

\vspace{2mm}

{\bf Lemma 3. (\cite{8})} Suppose $r>0,\,p>0$ and $\alpha>-1$. Then
there exists a constant $C>0$ such that

$$|f(z)|^{p}\leq\frac{C}{(1-|z|^{2})^{n+1+\alpha}}\int_{D(z,r)}|f(w)|^{p}d \upsilon_{\alpha}(w)$$
for all $f\in H(\textmd{B})$ and $z\in \textmd{B}.$

\vspace{2mm}

{\bf Lemma 4. (\cite{8})} Suppose $s>0$ and $t>-1$. Then
$$\int_{\textmd{S}}\frac{d \sigma(\xi)}{|1-\langle z,\xi\rangle|^{n+s}}\asymp (1-|z|^{2})^{-s}\eqno(6)$$
and
$$\int_{\textmd{B}}\frac{(1-|w|^{2})^{t}d \upsilon(w)}{|1-\langle z,w\rangle|^{n+1+t+s}}\asymp (1-|z|^{2})^{-s}\eqno(7)$$
as $|z|\rightarrow 1^{-}.$

It is known(see, e.g. \cite{8}) that for every $\delta$, there
exists a sampling sequence $\{a_{j}\}$ such that
$d(a_{j},a_{k})>\frac{\delta}{5}$ if $j\neq k$ and
$$\sum_{k=1}^{\infty} \chi_{D(a_{k},5\delta)}(z)\leq C. \eqno(8)$$

\vspace{2mm} {\bf Lemma 5. } Let $\sigma>1,\,t>0,\,\xi\in
\textmd{S},\,\widetilde{\Gamma}_{\sigma}(\xi)=\{z:\,|1-\overline{\xi}z|<\sigma(1-|z|)^{\frac{1}{n}}\}$.

Then there exist $\widetilde{\sigma}(\sigma,\,t)>1$ such that
$D(z,\,t)\subset \widetilde{\Gamma}_{\widetilde{\sigma}}(\xi)$ for
all $z\in \Gamma_{\sigma}(\xi).$ \vspace{2mm}

{\bf Proof of Lemma 5 } Let $w\in D(z,t),\,z\in
\Gamma_{\sigma}(\xi).$ We will show that $w\in
\widetilde{\Gamma}_{\widetilde{\sigma}}(\xi)$ for some
$\widetilde{\sigma}>1.$ Since $z\in \Gamma_{\sigma}(\xi),$ then
$|1-\overline{\xi}z|<\sigma(1-|z|),$ hence

\begin{eqnarray*}
  |1-\langle \xi, w\rangle|&\leq &  |1-\langle \xi, z\rangle|+  |\langle \overline{\xi},z\rangle-\langle \overline{\xi},w\rangle|\\
  &\leq &\sigma(1-|z|)+|z-w|\\
  &\leq &\sigma(1-|w|)+(\sigma+1)|z-w|.
\end{eqnarray*}

We will show $|z-w|\leq \sigma_{1}(1-|w|)^{\frac{1}{2}}$ for some
$\sigma_{1}>1.$ This is enough since $w\in D(z,t)$ is the same to
$z\in D(w,t)$ we have by exercise 1.1 from \cite{8}:
$$\frac{|P_{w}(z)-c|^{2}}{R^{2}\sigma_{1}^{2}}+\frac{|Q_{z}(w)|^{2}}{R^{2}\sigma_{1}}<1.\eqno(9)$$
where \begin{eqnarray*}&&R=\tanh(t),\,c=\frac{(1-R^{2})w}{1-R^{2}|w|^{2}},\,\sigma_{1}=\frac{1-|w|^{2}}{1-R^{2}|w|^{2}},\\
&&P_{w}(z)=\frac{\langle z,w\rangle
w}{|w|^{2}},\,Q_{w}(z)=z-\frac{\langle z,w\rangle w}{|w|^{2}}.
\end{eqnarray*}

Hence $$|z-w|\leq C_{1}\bigg(|z-P_{w}|+\bigg|\frac{\langle
z,w\rangle w}{|w|^{2}}-c\bigg|+|c-w|\bigg);\eqno(10)$$

$$|c-w|\leq |w|\bigg(1-\frac{1-R^{2}}{1-R^{2}|w|^{2}}\bigg)\leq \frac{C_{2}}{1-R^{2}}(1-|w|)=S_{2}.\eqno(11)$$

It is enough to show
$$\bigg|z-\frac{\langle z,w\rangle w}{|w|^{2}}\bigg|+\bigg|\frac{\langle z,w\rangle w}{|w|^{2}}-c\bigg|\leq R \bigg(\frac{1-|w|^{2}}{1-R^{2}|w|^{2}}\bigg)^{\frac{1}{2}}\widetilde{c}(R).\eqno(12)$$

Note $$|z-w|^{2}\leq C_{3} \bigg(|c-w|^{2}+\bigg|z-\frac{\langle
z,w\rangle w}{|w|^{2}}\bigg|^{2}+\bigg|\frac{\langle z,w\rangle
w}{|w|^{2}}-c\bigg|^{2}\bigg).\eqno(13)$$

Hence from (9)

\begin{eqnarray*}|z-w|^{2}&\leq& C_{4}(|c-w|^{2}+R^{2}\sigma_{1})\\
&\leq& C_{4}\bigg(S_{2}(|w|,R)^{2}+\frac{R^{2}(1-|w|^{2})}{1-R^{2}|w|^{2}}\bigg)\\
&\leq& C_{4}\bigg(S_{2}+\frac{R(1-|w|^{2})^{\frac{1}{2}}}{1-R^{2}|w|^{2}}\bigg)^{2}.
\end{eqnarray*}

Hence $|z-w|\leq  C_{5}(1-|w|)^{\frac{1}{2}}.$  Hence we complete
the proof of Lemma 5.

\vspace{2mm} {\bf Remark 2.} For $n=1$, Lemma 5 was proved by Wu's
paper in \cite{7}. \vspace{2mm}

{\bf Lemma 6. } Let $\mu$ be a positive Borel measure in
$\textmd{B}$. Let $D(w,\,t)\subset
\widetilde{\Gamma}_{\sigma}(\xi),\,w\in \Gamma_{\tau}(\xi),$ where
$\tau, \sigma,\,\xi,\,t$ are from Lemma 5. Let $f\in H(\textmd{B})$.
Then

$$\int_{\Gamma_{\tau}(\xi)}\frac{|f(z)|d\mu(z)}{(1-|z|)^{n}}\leq C\int_{\widetilde{\Gamma}_{\sigma}(\xi)}|f(z)|\int_{D(z,\,t)}d\mu(w)\frac{1}{(1-|z|)^{2n+1}}d\upsilon(z).$$

\vspace{2mm} {\bf Remark 3.} For $n=1$, Lemma 6 was provided in Wu's
paper in \cite{7}.

{\bf Remark 4.} The careful inspection of proof of Lemma 5 shows
that we can find $\delta,\,\delta>0$ such that $D(a,\,\delta)\subset
\widetilde{\Gamma}_{\widetilde{\sigma}}(\xi),$ if $a\in
\Gamma_{\sigma}(\xi),\,\sigma>1,\,$ for some fixed
$\widetilde{\sigma},\, \widetilde{\sigma}>1$.

\vspace{2mm}

{\bf Proof of Lemma 6} Since
$\chi_{D(z,\,t)}(w)=\chi_{D(w,\,t)}(z),\,z,w\in \textmd{B},\,t>0.$
Using Lemma 3, Lemma 5 and Fubini theorem, we have

\begin{eqnarray*}
\int_{\Gamma_{\tau}(\xi)}\frac{|f(z)|d\mu(z)}{(1-|z|)^{n}}
&\leq&C\int_{\Gamma_{\tau}(\xi)}\frac{1}{(1-|z|)^{n+1}}\bigg(\int_{D(z,\,t)}|f(w)|d\upsilon(w)\bigg)\frac{d\mu(z)}{(1-|z|)^{n}}\\
&\lesssim & \int_{\Gamma_{\tau}(\xi)} \frac{1}{(1-|z|)^{2n+1}}\int_{\textmd{B}}|f(w)|\chi_{D(w,\,t)}(z)d\upsilon(w)d\mu(z)\\
&\lesssim & \int_{ \widetilde{\Gamma}_{\sigma}(\xi)}\int_{\Gamma_{\tau}(\xi)} \frac{1}{(1-|z|)^{2n+1}}|f(w)|\chi_{D(w,\,t)}(z)d\upsilon(w)d\mu(z)\\
&\lesssim & \int_{ \widetilde{\Gamma}_{\sigma}(\xi)}\frac{1}{(1-|w|)^{2n+1}}|f(w)|\bigg(\int_{D(w,\,t)}d\mu(z)\bigg)d\upsilon(w).
\end{eqnarray*}

 Hence we complete the proof of Lemma 6.

\vspace{5mm}

We denote by $(A_{\alpha}^{p})_{1}$ the space of all holomorphic
functions in the unit ball  such that

$$\|f\|_{(A_{\alpha}^{p})_{1}}^{p}=\int_{\textmd{S}}\int_{\widetilde{\Gamma}_{t}(\xi)}\frac{|f(z)|^{p}(1-|z|)^{\alpha}}{(1-|z|)^{n}}
d\upsilon(z)d\sigma(\xi)<\infty,\,n\in\mathbb{N},\,0<p<\infty,\,\alpha>-1,$$
where $$\widetilde{\Gamma}_{t}(\xi)=\{z\in
\textmd{B}:\,|1-\xi\overline{z}|<t(1-|z|)^{\frac{1}{n}},\,t>1\}$$
 enlarged approach region.

Note for $n=1$, $(A_{\alpha}^{p})_{1}=A_{\alpha}^{p}.$ Since for
every $t>0,\,0<p<\infty,\,\alpha>-1,\,n\geq1,\,n\in\mathbb{N},$
$$\|f\|_{A_{\alpha}^{p}}^{p}\asymp\int_{\textmd{S}}\int_{\Gamma_{t}(\xi)}\frac{|f(z)|^{p}(1-|z|)^{\alpha}}{(1-|z|)^{n}}d\upsilon(z)d\sigma(\xi)$$
where $$\Gamma_{t}(\xi)=\{z\in
\textmd{B}:\,|1-\xi\overline{z}|<t(1-|z|)\},\,t>1,$$ we note that

\begin{eqnarray*}
\int_{\widetilde{\Gamma}_{t}(\xi)}|f(z)|^{p}d\upsilon_{\alpha}(z)
&\leq&\int_{\textmd{B}}|f(z)|^{p}d\upsilon_{\alpha}(z)\\
&\asymp& \int_{\textmd{S}}\int_{\Gamma_{t}(\xi)}\frac{|f(z)|^{p}d\upsilon_{\alpha}(z)}{(1-|z|)^{n}}d\sigma(\xi)\\
&\leq&\int_{\textmd{S}}\int_{\widetilde{\Gamma}_{t}(\xi)}\frac{|f(z)|^{p}d\upsilon_{\alpha}(z)}{(1-|z|)^{n}}d\sigma(\xi),
\end{eqnarray*}
where $0<p<\infty,\,\alpha>-1,\,f\in
(A_{\alpha}^{p})_{1}(\textmd{B}).$ We will use this observation in
the proof of Theorem\,1. \vspace{0.2cm}

We will also need Khinchine type estimate.

{\bf Remark 5.} Let us remind classical Khinchine's inequality. Let
$t\in [0,1),$
$$r_{0}(s)=\left \{\begin{array}{lll}&&\,\,\,\, 1,\,\,\,\,\,\,\,0\leq s-[s]<\frac{1}{2},\\
&&\\
&&-1,\,\,\,\,\,\,\,\frac{1}{2}\leq s-[s]<1.
\end{array}\right.$$

$$r_{j}(t)=r_{0}(2^{j}t),\,j=1,2,\ldots.\,\,\,\,\hbox{Classical\,\, Khinchine's\,\, inequality\,\,says: }~~~~~~~~~~~~~~~~~~~~~~~~~~~~$$
$$C_{p}\bigg(\sum_{j=0}^{N}|a_{j}|^{2}\bigg)^{\frac{p}{2}}\leq \int_{0}^{1}|\sum_{j=0}^{N}a_{j}r_{j}(t)|^{p}dt\leq \frac{1}{C_{p}}
\bigg(\sum_{j=0}^{N}|a_{j}|^{2}\bigg)^{\frac{p}{2}},\, 0<p<\infty,$$
Where $a_j$ are real numbers.

\vspace{0.2cm} We will also use the following assertion.

{\bf Remark 6.} Note that for any $\{z_{j}\}$-$r$-lattice in
$\textmd{B}$ and any large enough $m$ and for
$$f_{j}(z)=\frac{(1-|z_{j}|)^{m}}{(1-\overline{z_{j}}\,z)^{m+\frac{\alpha+n+1}{p}}},\,\,\,
\,\,\,j=1,2,\ldots,\,z\in \textmd{B}$$
for sufficient small $r>0,$
$\widetilde{f}(z)=\sum_{j}a_{j}f_{j}(z)$ is in
$A_{\alpha}^{p}(\textmd{B})$ by Theorem 2.30 of \cite{8}
 for any $a_{j}\in l^{p}$ and $\|\widetilde{f}\|_{A_{\alpha}^{p}}\leq C \|a_{j}\|_{l^{p}}.$

\section{Main results }

The goal of this section is to prove several direct generalizations
of recent results of Z. Wu from \cite{7} on the action of Luzin area
operator in Bergman classes in the unit disk. We consider such an
area operator based on ordinary Koranyi approach region and expanded
admissible approach region that coincide with each other in case of
unit disk and study it is action on
 Bergman classes in the unit ball.

\vspace{0.2cm}

{\bf Theorem 1. } Let $\mu$ be a positive Borel measure in
$\textmd{B}$. Let
 $f\in (A_{\alpha}^{p})_{1},\,p\leq q <\infty,\,0<p\leq1,\,\alpha>-1.$ Then
 $$\bigg\{\int_{\textmd{S}}\bigg(\int_{\Gamma_{t}(\xi)}\frac{|f(z)|d\mu(z)}{(1-|z|)^{n}}\bigg)^{q}
 d\sigma(\xi)\bigg\}^{\frac{1}{q}}\leq C\|f\|_{(A_{\alpha}^{p})_{1}},\,\,\hbox{for\,\,some\,\,}\,t>1$$
 if for some $\delta>0,$ $\int_{D(z,\delta)}d\mu(w)\leq C(1-|z|)^{\gamma},$ where $\gamma=\frac{\alpha+n+1}{p}+n-\frac{n}{q}.$
 The reverse assertion is true for $A_{\alpha}^{p}$ class and enlarged approach region.
\vspace{0.2cm}

 {\bf Remark 7.}
  For $n=1$, obviously enlarged approach region coinciding with ordinary approach
region, moreover $A_{\alpha}^{p}=(A_{\alpha}^{p})_{1}$, Theorem 1
was obtained by  Z. Wu in \cite{7}. \vspace{2mm}

{\bf Proof of Theorem 1 } First note that from Theorem 2.25 of
\cite{8} it follows that `` for some $\delta >0$ " in formulation of
Theorem 1  can be replaced by `` for any $\delta>0$ ".

 First we prove the second part. Let us note that the following estimate is true by Lemma 4.

$$\|f_{a}\|_{A_{\alpha}^{p}}\leq C,\,\,\, \hbox{if}\,\,\, f_{a}(z)=\frac{(1-|a|)^{m}}{(1-\overline{a}z)^{m+\frac{\alpha+n+1}{p}}},$$
where $a\in \textmd{B}$ is fixed, $m>0.$ Hence
$$G=\bigg\{\int_{\textmd{S}}\bigg(\int_{\widetilde{\Gamma}_{t}(\xi)}\frac{|f_{a}(z)|d\mu(z)}{(1-|z|)^{n}}\bigg)^{q}
 d\sigma(\xi)\bigg\}^{\frac{1}{q}}\leq  C.$$

 We now estimate $G=G(f_{a},\mu)$  from below. For that reason we use the following estimates.

 First if $z\in D(a,\delta)$, then by Lemma 2, we easily get $|f_{a}(z)|\asymp (1-|a|)^{-\frac{\alpha+n+1}{p}}$. On the other hand,
 the accurate inspection of proof of Lemma 5  (see Remark 4)
shows that we can find $\delta,\,\delta>0$ such that
$D(a,\,\delta)\subset \widetilde{\Gamma}_{\widetilde{\sigma}}(\xi),$
if $a\in \Gamma_{\sigma}(\xi),\,\sigma>1,\,$ for some fixed
$\widetilde{\sigma},\, \widetilde{\sigma}>1$. Using  all that we can
estimate $G(f_{a},\mu)$ from below to obtain the estimate we need.
We have

\begin{eqnarray*}
G(f_{a},\mu)^{q}&\geq& \int_{I_{\sigma}(a)}\bigg(\int_{D(a,\delta)}\frac{|f_{a}(z)|d\mu(z)}{(1-|z|)^{n}}\bigg)^{q}
 d\sigma(\xi)\\
&\geq& C (1-|a|)^{-\frac{\alpha+n+1}{p}q}\cdot (\mu(D(a,\delta)))^{q}\\
&&\cdot \frac{1}{(1-|a|)^{nq}}\cdot\int_{I_{\sigma}(a)} d\sigma(\xi),
\end{eqnarray*}

where $I_{\sigma}(a)=\{\xi\in \textmd{S}:\,a\in \Gamma_{\sigma}(\xi)
\},\,\sigma>0,$  $a\in \textmd{B}$. As it was noted in \cite{8}, the
following estimate holds.
$$|I_{\sigma}(a)|=\int_{I_{\sigma}(a)}d\xi=\int_{\textmd{S}}\chi_{\Gamma_{\sigma}(\xi)(a)}d\xi\asymp (1-|a|)^{n}. $$

Hence we have that $$G(f_{a},\mu)^{q}\geq
C(\mu(D(a,\delta)))^{q}(1-|a|)^{-\frac{\alpha+n+1}{p}q-nq+n}.$$

So finally we have

$$ \mu(D(a,\delta))\leq C(1-|a|)^{\frac{\alpha+n+1}{p}+n-\frac{n}{q}},$$

for some $\delta>0$ and all $a\in \textmd{B}.$

For $q=\infty,$ we have

\begin{eqnarray*}
\frac{\mu(D(a,\delta))}{(1-|a|)^{\frac{\alpha+n+1}{p}+n}}
&\leq& C\int_{\widetilde{\Gamma}_{\widetilde{\sigma}}(\xi)}\frac{|f_{a}(z)|d\mu(z)}{(1-|z|)^{n}}\\
&\leq& C\bigg\|\int_{\widetilde{\Gamma}_{\widetilde{\sigma}}(\xi)}\frac{|f_{a}(z)|d\mu(z)}{(1-|z|)^{n}}\bigg\|_{L^{\infty}}\\
&\leq& C  \|f_{a}\|_{A_{\alpha}^{p}}.
\end{eqnarray*}

The condition we obtained above on measure is also sufficient. To
show that we will need new estimates. For that reason we choose
$\{z_{j}\}$ lattice ($\widetilde{\delta}$-lattice) in $\textmd{B}$
with $\widetilde{\delta}$ less than $\delta,$
$$\mu(D(z_{j},\widetilde{\delta}))\leq C(1-|z_{j}|)^{\frac{\alpha+n+1}{p}+n-\frac{n}{q}},$$

Hence \begin{eqnarray*}
S_{p}(f,\mu)&=&\bigg(\int_{\Gamma_{t}(\xi)}\frac{|f(z)|d\mu(z)}{(1-|z|)^{n}}\bigg)^{p}\\
&\leq&  \bigg\{C\sum_{j,\,D(z_{j},\widetilde{\delta})\bigcap\Gamma_{t}(\xi)\neq \varnothing}\bigg(\sup_{w\in D(z_{j},\widetilde{\delta})}|f(w)|\bigg)\cdot \frac{\mu(D(z_{j},\widetilde{\delta}))}{(1-|z_{j}|)^{n}}\bigg\}^{p}\\
&\leq& C
\sum_{j,\,D(z_{j},\widetilde{\delta})\bigcap\Gamma_{t}(\xi)\neq
\varnothing}\bigg(\sup_{w\in
D(z_{j},\widetilde{\delta})}|f(w)|^{p}\bigg)\cdot
(1-|z_{j}|)^{\alpha+n+1-n\frac{p}{q}}.
\end{eqnarray*}

Furthermore $D(w,\widetilde{\delta})\subset
D(z_{j},2\widetilde{\delta})$ if $w\in D(z_{j},\widetilde{\delta})$
by triangle inequality for $\beta$ metric $\beta=\beta(z,w)$. And by
Lemma 3 and Lemma 2
$$S_{p}(f,\mu)<C \sum_{j}\int_{D(z_{j},2\widetilde{\delta})}|f(w)|^{p}(1-|w|)^{-\frac{np}{q}}d\upsilon_{\alpha}(w).$$

Since by Lemma 3 and Lemma 2, for $w\in D(z_{j},\widetilde{\delta})$
\begin{eqnarray*}
|f(w)|^{p}
&\leq& \frac{C}{(1-|w|)^{n+1+\alpha}}\int_{D(w,\widetilde{\delta})}|f(w)|^{p}d\upsilon_{\alpha}(w)\\
&\leq& \frac{C}{(1-|z_{j}|)^{n+1+\alpha}}\int_{D(z_{j},2\widetilde{\delta})}|f(w)|^{p}d\upsilon_{\alpha}(w).
\end{eqnarray*}

Since for $D(z_{j},\widetilde{\delta})\bigcap\Gamma_{t}(\xi)\neq
\varnothing$ and for $z\in
D(z_{j},\widetilde{\delta})\bigcap\Gamma_{t}(\xi),$ we have
$D(z_{j},2\widetilde{\delta})\subset D(z,3\widetilde{\delta})\subset
\widetilde{\Gamma}_{\widetilde{\sigma}}(\xi), \widetilde{\delta}>0,$
for some $\widetilde{\sigma}>0$, by Lemma 5 then using (8) we will
have

$$S_{p}(f,\mu)\leq C \int_{\widetilde{\Gamma}_{\widetilde{\sigma}}(\xi)}|f(z)|^{p}(1-|z|)^{-\frac{np}{q}}d\upsilon_{\alpha}(z)=M(f,p,q,\alpha)$$

for some $\sigma>0,$

$$(S_{p}(f,\mu))^{\frac{q}{p}}\leq C (M(f,p,q,\alpha))^{\frac{q}{p}}.$$

By H\"{o}lder inequality and the observation at the end of previous
section

$$(M(f,p,q,\alpha)))^{\frac{q}{p}}\leq \|f\|_{(A_{\alpha}^{p})_{1}}^{q-p}\int_{\widetilde{\Gamma}_{\widetilde{\sigma}}(\xi)}\frac{|f(z)|^{p}d\upsilon_{\alpha}(z)}{(1-|z|)^{n}}.$$

Hence integrating both sides by sphere $\textmd{S}$ we have finally
what we need.

$$\|(S_{p}(f,\mu))^{\frac{q}{p}}\|_{L^{1}(S)}\lesssim C\|f\|_{(A_{\alpha}^{p})_{1}}\asymp \int_{\textmd{S}} \int_{\widetilde{\Gamma}_{\widetilde{\sigma}}(\xi)}\frac{|f(z)|^{p}d\upsilon_{\alpha}(z)}{(1-|z|)^{n}}.$$

The proof of Theorem 1 is complete.

\vspace{5mm}

The following result follows from Lemma 5 and also can be found in
Wu's paper for $n=1$(see Theorem 3). \vspace{2mm}

{\bf Proposition 2.} Let $\alpha>-1,\,q \in (1,\infty),\, \mu$ is a
positive Borel measure in $\textmd{B},\,\delta>0,$

\begin{eqnarray*}&&\|A_{\mu}^{\delta}(f)\|_{L^{q}}=\bigg\{\int_{\textmd{S}} \bigg(\int_{\Gamma_{\delta}(\xi)}\frac{|f(z)|}{(1-|z|)^{n}}d\mu(z)\bigg)^{q}d\sigma(\xi)\bigg\}^{\frac{1}{q}}\leq C\sup_{z\in \textmd{B}}|f(z)|.
\end{eqnarray*}

Then
$$\bigg\|\int_{\widetilde{\Gamma}_{\tau}(\xi)}\bigg(\int_{D(z,\widetilde{\delta})}d\mu(w)\bigg)\frac{d\upsilon(z)}{(1-|z|)^{2n+1}}\bigg\|_{L^{q}(\textmd{S})}<\infty$$

for some $\widetilde{\delta}>0,\,\tau>0.$

{\bf Proof } By Lemma 5, for any $t$-lattice $\{z_{j}\}$ in
$\textmd{B}$, we have $D(z_{j},\,t)\subset
\widetilde{\Gamma}_{\tau}(\xi)$ or all $z_{j}\in
\Gamma_{\delta}(\xi)$, for some $\tau>0,\,\xi\in
\textmd{S},\,\delta>1,\,t>0.$ Hence since
$\|A_{\mu}^{\delta}(1)\|_{L^{q}}<\infty,$ we get what we need.
Indeed by (1),(2), Lemma 1 and Lemma 2 we have

\begin{eqnarray*}
&&\int_{\widetilde{\Gamma}_{\tau}(\xi)}\bigg(\int_{D(z,\widetilde{\delta})}d\mu(w)\bigg)\frac{d\upsilon(z)}{(1-|z|)^{2n+1}}\\
&\leq& C \sum_{j,\,D(z_{j},\widetilde{\delta})\bigcap\widetilde{\Gamma}_{\tau}(\xi)}\int_{D(z_{j},\,\widetilde{\delta})}\bigg(\int_{D(z_{j},2\widetilde{\delta})}d\mu(w)\bigg)
\frac{d\upsilon(z)}{(1-|z_{j}|)^{2n+1}}\\
&\leq& C \sum_{j}\bigg(\int_{D(z_{j},2\widetilde{\delta})}d\mu(w)\bigg)\frac{|D(z_{j},\,\widetilde{\delta})|}{(1-|z_{j}|)^{2n+1}}\\
&\leq& C  \sum_{j}\bigg(\int_{D(z_{j},2\widetilde{\delta})}d\mu(w)\bigg)\frac{1}{(1-|z_{j}|)^{n}}\leq C A_{\mu}^{\delta}(1).
\end{eqnarray*}
where at the last step we used arguments provided in proof of second
part of of theorem 1  which were based on estimate (8).

The proof of Proposition 2 is complete. \vspace{0.1cm}

The following result was obtained by Wu for $n=1$. We follow Wu's
ideas expanding his arguments to unit ball in $\mathbb{C}^n$.
\vspace{0.1cm}

{\bf Theorem 3.} Suppose $\alpha>-1,\,1<q<p\leq\infty,\,t>1,$ $\mu$
be a nonnegative Borel measure in $\textmd{B}$. Then the following
is true.
$$\bigg\{\int_{\textmd{S}}\bigg(\int_{\Gamma_{t}(\xi)}\frac{|f(z)|d\mu(z)}{(1-|z|)^{n}}\bigg)^{q}
 d\sigma(\xi)\bigg\}^{\frac{1}{q}}\leq C\|f\|_{(A_{\alpha}^{p})_{1}},$$
 if
 $$K_{p,q}=\int_{\widetilde{\Gamma}_{\widetilde{\sigma}}(\xi)}\bigg(\int_{D(z,\delta)}
 d\mu(w)\bigg)^{\frac{p}{p-1}}\frac{d\upsilon(z)}{(1-|z|)^{n+(n+1)\frac{p}{p-1}}} \in L^{\frac{(p-1)q}{p-q}}(\textmd{S})$$
 for any $\widetilde{\sigma}>1,\,\delta>0.$

 {\bf Proof }  Let us consider first $p=\infty $ case. The proof follows directly from Lemma 6 that we proved above. Let $f\in H(\textmd{B})$. Then  we easily have by Lemma 6
 \begin{eqnarray*}
M_{q}(f,\mu)&=&\int_{\textmd{S}}\bigg(\int_{\Gamma_{t}(\xi)}\frac{|f(z)|d\mu(z)}{(1-|z|)^{n}}\bigg)^{q}
 d\sigma(\xi)\\
 &\leq& C\int_{\textmd{S}}\bigg(\int_{\widetilde{\Gamma}_{\widetilde{\sigma}}(\xi)}\frac{|\widetilde{f}(z)|d\upsilon(z)}{(1-|z|)^{2n+1}}\bigg)^{q} d\sigma(\xi),
\end{eqnarray*}
where $$|\widetilde{f}(z)|=\bigg(\int_{D(z,\delta)}
 d\mu(w)\bigg)\cdot|f(z)|.$$

 Hence $$M_{q}(f,\mu)\leq \bigg(\sup_{z\in \textmd{B}}|f(z)|\bigg)\cdot K(\mu). $$

 This is enough since $K_{\infty,q}(\mu)=K(\mu).$

 Let $q>1,\,1<p<\infty.$ Then again using Lemma 6 and  H\"{o}lder inequality, we have the following chain of estimates.
\begin{eqnarray*}
S(\mu,f,t)&=&\int_{\Gamma_{t}(\xi)}\frac{|f(z)|d\mu(z)}{(1-|z|)^{n}}\\
 &\leq& C \int_{\widetilde{\Gamma}_{\sigma}(\xi)}\frac{|f(z)|\bigg(\int_{D(z,\delta)}
 d\mu(w)\bigg)d\upsilon(z)}{(1-|z|)^{2n+1}}\\
  &\lesssim & C \bigg(\int_{\widetilde{\Gamma}_{\sigma}(\xi)}\frac{|f(z)|^{p}d\upsilon(z)}{(1-|z|)^{n}}\bigg)^{\frac{1}{p}}\\
  &&\cdot
  \bigg\{\int_{\widetilde{\Gamma}_{\sigma}(\xi)}\frac{1}{(1-|z|)^{q}}\bigg(\int_{D(z,\delta)}
 d\mu(w)\bigg)^{p^{'}}\frac{d\upsilon(z)}{(1-|z|)^{n}}\bigg\}^{\frac{1}{p^{'}}}.
\end{eqnarray*}

Hence again using H\"{o}lder inequality,
$q=(n+1)p^{\prime},\,\frac{1}{p}+\frac{1}{p^{\prime}}=1,$
$$\int_{\textmd{S}}(S(\mu,f,t))^{q}d\sigma(\xi)\leq C \|f\|^{q}_{(A_{\alpha}^{p})_{1}}.$$

 The proof of Theorem 3 is complete.

The following result was also obtained by Wu for $n=1.$  We follow
again Wu's ideas to expand his one-dimensional arguments to the unit
ball case in $\mathbb{C}^n$.

{\bf Theorem 4.} Let $\alpha>-1,\,q<p\leq\infty,\,0<q<1.$ Let $\mu$
be a nonnegative Borel measure in $\textmd{B}$. Then the following
statements are true.

1) For any fixed $\sigma>0,$
$$\bigg\{\int_{\textmd{S}}\bigg(\int_{\Gamma_{\sigma}(\xi)}\frac{|f(z)|d\mu(z)}{(1-|z|)^{n}}\bigg)^{q}
 d\sigma(\xi)\bigg\}^{\frac{1}{q}}\leq C\|f\|_{A_{\alpha}^{p}},\eqno(M_{1})$$
 if $$\int_{\textmd{S}}\bigg(\sum_{j,\,z_{j}\in \widetilde{\Gamma}_{\tau}(\xi)}|a_{j}|\int_{D(z_{j},\delta)}d\mu(w)\frac{1}{(1-|z_{j}|)^{s}}
\bigg)^{q}d\sigma(\xi)\leq \|\{a_{j}\}\|_{l^{p}}^{q}\eqno(M_{2})$$
for any $\tau>0,\,\delta-$lattice $\{z_{j}\}$,
$s=\frac{\alpha+n+1}{p}+n.$

2) If ($M_{1}$) holds  for enlarged Koranyi region, then ($M_{2}$)
holds  for some $\tau>0$ and ordinary  admissible Koranyi region.

\vspace{0.02cm}

{\bf Remark 8.} For $n=1,$  Theorem 4 was proved by  Wu's theorem in
\cite{7}. \vspace{0.02cm}

{\bf Proof of Theorem 4 }  Let first condition  ($M_{2}$) holds. As
in proof of Theorem 1 (arguments in second part of proof based on
estimate (8)), $D(z_{j},\delta)\bigcap \Gamma_{\sigma}(\xi)\neq
\emptyset$ implies $z_{j}\in \widetilde{\Gamma}_{\tau}(\xi)$ for
some $\tau>1$ (Note this for n=1 was proved in \cite{7}).  Hence for
any function $f$, $f\in H(\textmd{B})$, we have the following chain
of estimates, let $D_{j}=D(z_{j},\delta),$
\begin{eqnarray*}
\int_{\Gamma_{\sigma}(\xi)}\frac{|f(z)|d\mu(z)}{(1-|z|)^{n}}
&\lesssim & C \sum_{j,\,D_{j}\bigcap \Gamma_{\sigma}(\xi)\neq \emptyset}\int_{D_{j}} \frac{|f(z)|d\mu(z)}{(1-|z|)^{n}}\\
 &\lesssim & C\sum_{j} \bigg(\sup_{z\in D_{j}} |f(z)|\bigg) \frac{\mu(D_{j})}{(1-|z_{j}|)^{n}}\\
  &\lesssim & C\sum_{j} \bigg(\sup_{z\in D_{j} }|f(z)|\bigg) (1-|z_{j}|)^{\frac{\alpha+n+1}{p}}\\
  && \cdot
  \frac{1}{(1-|z_{j}|)^{(\alpha+n+1)(\frac{1}{p}-1)+n}}\bigg(\int_{D(z_{j},\delta)}d\mu(w)\bigg)\\
  &&\cdot\frac{1}{(1-|z_{j}|)^{\alpha+n+1}}\\
   &\lesssim & C\sum_{j,\,z_{j}\in \widetilde{\Gamma}_{\tau}(\xi)} \sup_{z\in D_{j}}  |f(z)| (1-|z_{j}|)^{\frac{\alpha+n+1}{p}}\\
   &&\cdot \bigg(\int_{D(z_{j},\delta)}d\mu(w)\bigg)\cdot\frac{1}{(1-|z_{j}|)^{\frac{\alpha+n+1}{p}+n}}.
\end{eqnarray*}

Hence we only need to show that
$$(1-|z_{j}|)^{\frac{\alpha+n+1}{p}}\sup_{z\in D_{j}}  |f(z)|$$ is
in $l^{p}$ and it is norm dominated by $ C\|f\|_{A_{\alpha}^{p}}$.
This is obvious for $A_{\alpha}^{\infty}=H^{\infty},$ where
$H^\infty$ is the class of all bounded analytic functions in
$\textmd{B}$.  For $p<\infty,$ we have using Lemma 1-3

$$\sum_{j}(1-|z_{j}|)^{\alpha+n+1}\bigg(\sup_{z\in D_{j}}|f(z)|^{p} \bigg)\leq C\sum_{j}\int_{D(z_{j},2\delta)}|f(w)|^{p}d\upsilon_{\alpha}(w)\lesssim C \|f\|_{A_{\alpha}^{p}}^{p}.$$

This is what we need.

Let us show the reverse assertion in Theorem 4. We will need several
auxiliary assertions for that. First a theorem from \cite{3}, on
Khinchine's inequality for $K-$quasinorm $\|\cdot\|$ of $K-$qausi
Banach space $X$, then we will need the atomic decomposition of
Bergman classes $A_{\alpha}^{p}$
 which can be found in \cite{8} and we will need also Lemma 5.

Recall that a $K-$quasinorm $\|\cdot\|$ for quasi-Banach space $X$
is a function from $X$ to $[0,\infty)$ which has all the properties
of a norm except that the triangle inequality is replaced by
$\|x+y\|\leq K (\|x\|+\|y\|)$ for all $x,y\in X.$

{\bf Theorem C (see \cite{3})} Suppose $0<\tau<s<\infty, K>0.$ There
exist positive constants $C$ and $C$ depending only on $\tau,s,K$
such that for any quasi-Banach space $X$ with $K-$quasinorm
$\|\cdot\|$, any positive integer $n$ and any $x_{1},\ldots,x_{n}\in
X$, the following estimates hold.
\begin{eqnarray*}
C\bigg(\int_{0}^{1}\bigg\|\sum_{j}r_{j}(t)x_{j}\bigg\|^{\tau}dt\bigg)^{\frac{1}{\tau}}
 &\leq& \bigg(\int_{0}^{1}\bigg\|\sum_{j}r_{j}(t)x_{j}\bigg\|^{s}dt\bigg)^{\frac{1}{s}}\\
   &\leq& C\bigg(\int_{0}^{1}\bigg\|\sum_{j}r_{j}(t)x_{j}\bigg\|^{\tau}dt\bigg)^{\frac{1}{\tau}}.
\end{eqnarray*}

{\bf Theorem D (see \cite{8})} Suppose $p>0,\alpha>-1$ and
$l>n\max\{1,\frac{1}{p}\}+\frac{\alpha+1}{p}.$ Then there exists a
sequence $\{z_{k}\}$ in $\textmd{B}$ such that
$A_{\alpha}^{p}(\textmd{B})$ consists exactly of functions of the
form

$$f(z)=\sum_{k=1}^{\infty}a_{k}\frac{(1-|z_{k}|^{2})^{(lp-n-1-\alpha)/p}}{(1-\langle z,z_{k}\rangle)^{l}},\,z\in \textmd{B},\,\,
\|a_{k}\|_{l^{p}}\asymp \|f\|_{A_{\alpha}^{p}}.$$
where $a_{k}$ belongs to the sequence space $l^{p}$ and the series
converges in the norm
 topology of $A_{\alpha}^{p}(\textmd{B})$.

 \vspace{0.2cm}
 {\bf Remark 9.} The sequence ${z_{k}}$ from Theorem D is a $r$-lattice, for small enough $r>0$ (  see \cite{8}).

 Now we turn to the proof of Theorem 4, expanding Wu's arguments to the unit ball.

 Assume that $$\bigg\{\int_{\textmd{S}}\bigg(\int_{\widetilde{\Gamma}_{\sigma}(\xi)}\frac{|f(z)|d\mu(z)}{(1-|z|)^{n}}\bigg)^{q}
 d\sigma(\xi)\bigg\}^{\frac{1}{q}}\leq C\|f\|_{A_{\alpha}^{p}}\,\,\,\,\hbox{ for\,\, any }\,\,\sigma>0.\eqno(14)$$

 Let $D_{\rho}=\{z\in \textmd{B}:\,|z|<\rho\}$, $\rho\in(0,1)$. Let $f\in H(\textmd{B}).$ Put

$$ \|f\|_{\xi,\sigma,\rho}=\int_{\widetilde{\Gamma}_{\sigma}(\xi)\bigcap D_{\rho}}\frac{|f(z)|d\mu(z)}{(1-|z|)^{n}}.$$

Obviously for $q\in (0,1)$ by Theorem C we have for any $f_{j}\in
A_{\alpha}^{p}(\textmd{B})$

$$\bigg(\int_{0}^{1}\bigg\|\sum_{j}a_{j}r_{j}(t)f_{j}(z)\bigg\|_{\xi,\sigma,\rho}^{q}dt\bigg)^{\frac{1}{q}}
\asymp \int_{0}^{1}\bigg\|\sum_{j}a_{j}r_{j}(t)f_{j}(z)\bigg\|_{\xi,\sigma,\rho}dt$$

Using Fubini's theorem  and Khinchine's inequality (see also Remark
5) we have
$$\int_{0}^{1}\bigg|\sum_{j}a_{j}r_{j}(t)f_{j}(z)\bigg|dt\asymp \bigg(\sum_{j}|a_{j}|^{2}|f_{j}(z)|^{2}\bigg)^{\frac{1}{2}}.$$

$$\!\int_{\textmd{S}}\bigg(\int_{\widetilde{\Gamma}_{\sigma}(\xi)\bigcap D_{\rho}}\bigg(\sum_{j}|a_{j}|^{2}|f_{j}(z)|^{2}\bigg)^{\frac{1}{2}}
\frac{d\mu(z)}{(1-|z|)^{n}}\bigg)^{q}d\sigma(\xi)
\!\!\asymp\!\!\int_{\textmd{S}}\int_{0}^{1}\bigg\|\sum_{j}a_{j}r_{j}(t)f_{j}(z)\bigg\|_{\xi,\sigma,\rho}^{q}dt d\sigma(\xi).\eqno(15)$$

It is easy to note that constants in equivalence relation (15) do
not depend on $\rho$ so we can pass to limit $\rho\rightarrow1$.

Let $$\{a_{j}\}\in
l^{p};\,\,\,\,f_{j}(z)=\frac{(1-|z_{j}|)^{m}}{(1-\overline{z_{j}}\,z)^{m+\frac{\alpha+n+1}{p}}}\,\,\,
\hbox{with}\,\,\, m>n+1,\,j=1,2,\ldots,\,z\in \textmd{B}.$$

 Then by Theorem D and Remark 6
$$f_{t}(z)=\sum_{j}a_{j}r_{j}(t)f_{j}(z)\in A_{\alpha}^{p}(\textmd{B})\,\,\,\,\hbox{and}\,\,\,\,
 \|f_{t}\|_{A_{\alpha}^{p}}\leq C \|\{a_{j}\}\|_{l^{p}},\,t\in[0,1).$$

 Putting this $f_t$ into (14) and (15) and using the fact that

$$\bigg(\sum_{j}|a_{j}|^{2}|f_{j}(z)|^{2}\bigg)^{\frac{1}{2}}\geq |a_{j}||f_{j}(z)|\asymp \frac{|a_{j}|}{(1-|z_{j}|)^{m+\frac{\alpha+n+1}{p}}},\,
z_{j}\in D(z,\widetilde{\delta}),\,\widetilde{\delta}>0,$$

we have
$$ \int_{\textmd{S}}\bigg(\sum_{j,D_{j}\subset \widetilde{\Gamma}_{\sigma}(\xi)}\frac{|a_{j}|\mu(D_{j})}{(1-|z_{j}|)^{\frac{\alpha+n+1}{p}+n}}
\bigg)^{q}d\sigma(\xi)
 \leq C \|\{a_{j}\}\|^{q}_{l^{p}}.$$

Hence\,\, we\,\, will\,\, have
$$  \int_{\textmd{S}}\bigg\{\sum_{j,z_{j}\in \Gamma_{\tau}(\xi)}|a_{j}|\bigg(\int_{D(z_{j},\delta)}d\mu(w)\bigg)
\frac{1}{(1-|z_{j}|)^{\frac{\alpha+n+1}{p}+n}}\bigg\}^{q}d\sigma(\xi)
  \leq C \|\{a_{j}\}\|^{q}_{l^{p}}.$$

 Since as it was noted before if $z_{j}\in \Gamma_{\tau}(\xi)$ for some $\tau>0,$
 then $D(z_{j},\delta)\subset \widetilde{\Gamma}_{\sigma}(\xi)$.

  The proof of Theorem 4 is complete.

\vspace{0.2cm}

{\bf Remark 10.} It is easy to notice that the proof of the first
part of Theorem 4 is true for all $p>0$ and $ q>0$.

\end{document}